\documentclass[11pt]{article}
\usepackage{amsfonts, latexsym, amssymb,amsmath,amstext}
\usepackage[dvips]{graphics,color}
\usepackage[mathscr]{eucal}
\newenvironment{proof}{\par\noindent\emph{Proof.}\ }{\hfill$\Box$\par}
\newtheorem{theorem}{Theorem}[section]

\newtheorem{remark}{Remark}[section]

\newtheorem{proposition}{Proposition}[section]
\makeatletter
\@addtoreset{equation}{section}

\makeatother
\begin{document}
\title{Remarks Towards the Classification of $RS_4^2(3)$-Transformations and
Algebraic Solutions\\ of the Sixth Painlev\'e Equation}
\author{A.~V.~Kitaev
\thanks{E-mail: kitaev@pdmi.ras.ru}\\
Steklov Mathematical Institute, Fontanka 27, St.Petersburg, 191023, Russia\\
and\\
School of Mathematics and Statistics, University of Sydney,\\
Sydney, NSW 2006, Australia}
\date{February 21, 2005}
\maketitle
\begin{abstract}
We introduce a notion of the divisor type for rational functions 
and show that it can be effectively used for the classification of 
the deformations of dessins d'enfants related with the 
construction of the algebraic solutions of the sixth Painlev\'e 
equation via the method of $RS$-transformations.
\vspace{24pt}\vspace{24pt}\\
{\bf Short title}: {Classification of $RS_4^2(3)$-Transformations}
\end{abstract}
\newpage
\setcounter{page}2
\section{Introduction}
  \label{sec:intro}
Recently the author introduced a general method of $RS$-trans\-for\-ma\-tions~\cite{K2}  
for special functions of the isomonodromy type (SFITs)~\cite{K1}. This method applies to 
SFITs defining isomonodromy deformations of linear $n\times n$-mat\-rix ODEs of the first 
order with rational coefficients and with both regular and essential singular points. 

$RS$-Trans\-for\-ma\-tions are just a proper combination of rational transformations 
($R$-trans\-for\-ma\-tions) of the independent variable of the linear ODEs and 
Schlesinger transformations ($S$-trans\-for\-ma\-tion) of the dependent variable.
Solutions of many different and seemingly unrelated problems from various areas of the 
theory of functions get a unified and systematic approach in the framework of this method 
and can be reduced to the study, construction, and classification of different
$RS$-trans\-for\-ma\-tions for matrix linear ODEs.
This method, e.g., allows one to prove the duplication formula for the Gamma-function 
(and most probably the general multiplication formula for the multiple argument~\cite{AAR}),
build higher-order transformations for the Gauss hypergeometric function and reproduce the
Schwarz table for it \cite{AK1,HGBAA}, construct quadratic transformations for the 
Painlev\'e and classical transcendental functions \cite{K4,K5}, and provide a systematical
method for finding algebraic points at which transcendental SFITs attain algebraic 
values~\cite{AK2}. Without doubt, many other interesting problems can be approached via 
the method of $RS$-trans\-for\-ma\-tions. In this paper we apply this general method to 
the problem of construction and classification of algebraic solutions of the sixth 
Painlev\'e equation. 

Recently scanning the literature, I realized that, possibly, the first serious profound  
result concerning $RS$-trans\-for\-ma\-tions was obtained by F.~Klein~\cite{Klein84}, who 
proved that any scalar Fuchsian equation of the
second order with finite monodromy group is a "pull-back" ($R$-trans\-for\-ma\-tion) of 
the Euler hypergeometric equation. In this context instead of the $S$-trans\-for\-ma\-tions
the notion of "projective equivalence" is used. The latter is more restrictive than 
general $S$-trans\-for\-ma\-tions because in terms of the matrix ODEs it corresponds to 
triangular Schlesinger transformations, that finally results in a more restrictive special 
choice of the exponent differences (formal monodromy) of the hypergeometric equation, than
when more general $S$-trans\-for\-ma\-tions are allowed.

Klein's result immediately implies that any solution of the Garnier system and, in 
particular the sixth Painlev\'e equation that corresponds to a finite monodromy group of 
the associated Fuchsian equation, is algebraic. It is important to mention that the 
converse statement is not true.

In the context of the sixth Painlev\'e equation the first person who could, theoretically, 
apply the "pull-back ideology" was R.~Fuchs because it was he who found that the sixth 
Painlev\'e equation governs isomonodromy deformations of the certain scalar second order 
Fuchsian ODE and, moreover, received an informative letter from F.~Klein. He actually did 
it, in a study of algebraic solutions in the so-called Picard case of the sixth Painlev\'e 
equation~\cite{F1,F2}\footnote{These works were not known to me and,
possibly, to most modern researchers until very recently, when Yousuke Ohyama called our 
attention to them}.
Recently appeared a paper by Ch.~Doran~\cite{D} who formulated a more general scheme for 
construction of algebraic solutions of the sixth Painlev\'e equation from the pull-back 
point of view. A more detailed account of the last work is given in Introduction of 
\cite{HGBAA}. In the following two paragraphs we explain why the method of 
$RS$-trans\-for\-ma\-tions for construction of the algebraic solutions is more general than 
the pull-back back one.

For a given $R$-trans\-for\-ma\-tion one can normally associate a few different 
$RS$-trans\-for\-ma\-tions, due to the possibility of choosing different (not related by 
the contiguity transformations) initial hypergeometric equations which suffer this
$R$-trans\-for\-ma\-tion and, by further application of proper $S$-trans\-for\-ma\-tions,
are mapped into the Fuchsian ODE with four regular points. Each of these 
$RS$-trans\-for\-ma\-tions generate an algebraic solution of the sixth Painlev\'e equation,
which sometimes depends on a complex parameter. Thus we have a finite number of algebraic
solutions. The subset of these solutions belonging to different orbits of the group of
$RS$-trans\-for\-ma\-tions with $\deg R=1$, we call, {\it seed algebraic solutions}, and 
their generating $RS$-trans\-for\-ma\-tions - {\it seed $RS$-trans\-for\-ma\-tions}.
The seed algebraic solutions corresponding to the same rational covering 
($R$-transformation) are different, by  definition; however, the seed solutions associated 
with different rational coverings can coincide. Furthermore, the seed solutions, even 
corresponding to the same rational covering, can sometimes be related by some compositions 
of the quadratic transformations and/or B\"acklund transformations. Since the quadratic 
transformations are generated by the $RS$-trans\-for\-ma\-tions with $\deg R=2$, and one 
of the B\"acklund transformations has no realization as the Schlesinger transformation 
of the $2\times2$-matrix Fuchsian ODE; we call this special transformation the Okamoto 
transformation (see \cite{O} and Appendix~\cite{HGBAA,KK}).

We call attention of the reader that the possibility of the construction of different 
$RS$-trans\-for\-ma\-tions starting from the same rational covering mentioned in the 
previous paragraph is not considered by the successors of the "pull-back ideology"
because of the projective invariance property which assumes only one particular choice of 
the formal monodromy of the initial hypergeometric equation. Therefore, the "pull-back 
results" in many cases, namely in those ones where the property of projective equivalence 
can be changed to a less restrictive condition of the existence of $S$-trans\-for\-ma\-tion, 
can be extended or completed. We discuss this opportunity for construction of higher-order 
transformations of the Gauss hypergeometric functions in the Remarks 
in Sections~\ref{sec:D7} and~\ref{sec:D8}.

This paper is a continuation of author's previous work~\cite{HGBAA}. In \cite{HGBAA} 
we give a general definition of the one-dimensional deformations of dessins d'enfants and 
their relation to the algebraic solutions of the sixth Painlev\'e equation, construct by 
this method numerous examples of different algebraic solutions, and discuss different
features of this technique, e.g., a mechanism of appearance of genus-$1$ algebraic 
solutions.
Here we put this technique onto a systematic footing. A new idea we use here is symmetry 
preserving and symmetry braking deformations of the dessins d'enfants and their relation 
to uniqueness of the corresponding rational covering.

More precisely, in Section~\ref{sec:main} we introduce a notion of {\it divisor type} 
($D$-type) of a rational function and classify all $D$-types of the rational functions 
that generate algebraic solutions of the sixth Painlev\'e equation via the method of 
$RS$-transformations ($R_4(3)$-{\it functions}). We call the $D$-{\it series} the set of 
all $R_4(3)$-functions having the same $D$-type. Among these $D$-series there are two ones 
with finitely many, actually a few, members.
This fact is proved and the corresponding rational
functions are explicitly constructed in Sections~\ref{sec:D7} and~\ref{sec:D8}.
Each of the other $D$-series, corresponding to the $D$-types specified in the 
classification theorem of Section~\ref{sec:main}, are infinite.

It is worth noticing that modern personal computers (PC) allows one to construct all 
rational coverings that are presented here and in~\cite{HGBAA} without any advanced 
algorithms just by the natural method explained in Remark~$2.1$ of~\cite{HGBAA}. The time 
of calculation with MAPLE code on a relatively powerful PC does not exceed $1$ second for 
any of these functions. 
Of course, finding the concise parametrization requires much more 
additional time. It is interesting to note that in 1998-2000, when we used exactly
the same calculational scheme but on the Pentium 2 based PC with about 256 Mb RAM, we were
not able to construct many interesting functions, even some Belyi function of degree $8$,
see~\cite{AK1} we have found only numerically. 
This remark, however, does not mean that we do not need any advanced 
calculational algorithms; explicit construction of most of the rational coverings with 
the degree $>12$ still represent substantial difficulties.  
   
To each $R_4(3)$-function we also indicate the number of the seed $RS$-transformations and
present one algebraic solution whose construction does not require explicit form of the 
related Schlesinger transformation. It is exactly the "pull-back" solution to get explicitly
the other seed solutions one has to construct explicitly corresponding $S$-transformations.
This procedure is absolutely straightforward and does not require any advance computer
algorithms and we do not consider it here. Numerous examples of the complete constructions
of $RS$-transformations are given in \cite{AK2}.   

This paper is a far-going extension of the second part of my talk in Angers, where I have 
only explained some simplest ideas concerning the concept of deformations of the dessins 
d'enfants and announced the construction of the solution presented in Section~\ref{sec:D7}.

During the preparation of this paper there appeared two papers by P.~Boalch~\cite{Bo1,Bo2}, 
who is classifying algebraic solutions of the sixth Painlev\'e equation by developing the 
method (or, perhaps, more precisely to say following the trend)
suggested by B.~A.~Dubrovin and M.~Mazzocco~\cite{DM}. In view of author's conjecture,
that all algebraic solutions can be obtained from the seed algebraic solutions 
with the help of quadratic and B\"acklund transformations, it is important to mention
that all solutions which are specified in these works satisfy this hypothesis
(actually many of them are equivalent or related with the ones already published 
in \cite{AK2,HGBAA}). The other (most of them are already obtained) will appear in the 
papers devoted to systematic studies of the infinite $D$-series.\\
\noindent
{\bf Acknowledgment} The author is grateful to Michele Loday and Eric Delabaere,
the organizers of the conference in Angers, for the invitation and prompt resolution of
the organizational problems allowing him to participate in the conference,
Kazuo Okamoto for the invitation to the University of Tokyo, where this work was
finished, and hospitality, Yousuke Ohyama for the invitation to the University of Osaka,
where the final version of this work was presented, and for references~\cite{F1,F2},
Hidetaka Sakai for valuable discussions and various help during his stay in Tokyo,
and Philip Boalch for many discussions of different aspects concerning algebraic solutions.\\
The work is supported by JSPS grant-in-aide no.~$14204012$.
\section{$D$-Type of Rational Functions}
 \label{sec:main}
We begin this part of the lecture
with the canonical form of the sixth Painlev\'e equation, because
we are going to present a few algebraic solutions of this equation
in the explicit form.
\begin{eqnarray}
 \label{eq:P6}
\frac{d^2y}{dt^2}&=&\frac 12\left(\frac 1y+\frac 1{y-1}+\frac 1{y-t}\right)
\left(\frac{dy}{dt}\right)^2-\left(\frac 1t+\frac 1{t-1}+\frac 1{y-t}\right)
\frac{dy}{dt}\nonumber\\
&+&\frac{y(y-1)(y-t)}{t^2(t-1)^2}\left(\alpha_6+\beta_6\frac t{y^2}+
\gamma_6\frac{t-1}{(y-1)^2}+\delta_6\frac{t(t-1)}{(y-t)^2}\right),
\end{eqnarray}
where $\alpha_6,\,\beta_6,\,\gamma_6,\,\delta_6\in\mathbb C$ are
parameters. For a convenience of comparison of the results obtained here
with the ones from the other works we will use also parametrization of the
coefficients in terms of the formal monodromies $\hat\theta_k$:
$$
\alpha_6=\frac{(\hat\theta_\infty-1)^2}2,\quad
\beta_6=-\frac{{\hat\theta}_0^2}2,\quad\gamma_6=\frac{{\hat\theta}_1^2}2,
\quad\delta_6=\frac{1-{\hat\theta}_t^2}2.
$$
It is well known that any solution of this equation defines an
isomonodromy deformation of the $2\times2$ matrix Fuchsian ODE on
the Riemann sphere with four singular points.

As a first step in construction of algebraic solutions of
Equation~(\ref{eq:P6}) via the method of $RS$-transformations one
has to construct a proper rational covering of the Riemann sphere.
The corresponding rational function has four critical values.
Three of them are supposed to be placed at $0$, $1$, and $\infty$.
To specify proper rational functions we use the symbol of their {\it $R$-type},
which consists of three boxes. In these boxes we consecutively write
partitions of multiplicities of preimages of the points $0$, $1$, and $\infty$,
correspondingly. The fourth critical point has a standard partition of its
multiplicities $2+1+\ldots+1$ which is not indicated in the $R$-type.

According to \cite{HGBAA} the numbers in each box can be presented as a
union of two nonintersecting sets: the {\it apparent} set and {\it nonapparent} one.
The characteristic property of the apparent set is that g.c.d. of its members
is $\geq2$. It might be that nonapparent set has also nontrivial g.c.d.,
thus in general a presentation of the box as a union of the apparent and
nonapparent sets is not unique. Moreover, nonapparent set may contain a number
which is divisible by the g.c.d. of the apparent set. When the subdivision of
the boxes in the apparent and nonapparent sets is chosen we have an ordered triplet of
three integer numbers, $<m_0,m_1,m_\infty>$, the divisors of the apparent sets in the
corresponding boxes, which we call the {\it divisor type
{\rm(}$D$-type{\rm)} of the rational function}.
\begin{remark}
 \label{rem:convention} 
{\rm In our notation of the $D$-types we always assume that $m_0\leq m_1\leq m_\infty$,
clearly this can always be achieved by rearranging the points $0$, $1$, and $\infty$ by
a fractional-linear transformation. However, in the notation of $R$-types we do not follow 
this agreement, and in most cases we have $m_0>m_\infty>m_1$. Actually, we can speak of
the two types of numbering of the boxes in $R$-types: the natural one, i.e., according to
their position in the $R$-sym\-bol; and the $D$-con\-sis\-tent numbering, i.e., according to 
the rule: the larger g.c.d., the larger number of the box. In the statements we always
assume the $D$-con\-sis\-tent numbering in the proofs - the natural one.}
\end{remark}
Another important parameter of the proper rational functions is the total number
of members in all three nonapparent sets. In our case this number is $4$. We put
this number as the subscript in the notation of the $R$-type:
$R_4(\ldots|\ldots|\ldots)$ or in short $R_4(3)$. To simplify notation we omit the
subscript in situations where it cannot course a confusion.
\begin{theorem}
 \label{thm:main}
$R_4(3)$-Rational functions have one of the following eight $D$-types:\\
$<2,2,m>$, $<2,3,3>$, $<2,3,4>$, $<2,3,5>$, $<2,3,6>$, $<2,3,7>$,
$<2,3,8>$, $<2,4,4>$. Where $m-1\in\mathbb N$.
\end{theorem}
\begin{proof}
Let $n\geq2$ be the degree of some $R_4(3)$ function of $D$-type $<m_0,m_1,m_\infty>$.
Denote the sum of numbers in nonapparent sets in the consecutive boxes of a $R_4(3)$
function as $\sigma_0$, $\sigma_1$, and $\sigma_\infty$, respectively.

From the Riemann-Hurwitz formula, with a help of Proposition~2.1
of \cite{HGBAA}, one deduces the following "master" inequality,
\begin{equation}
 \label{ineq:master}
\frac{n-\sigma_0}{m_0}+\frac{n-\sigma_1}{m_1}+
\frac{n-\sigma_\infty}{m_\infty}\geq n-1.
\end{equation}
Clearly the numbers $\sigma_k$ satisfy one more inequality
\begin{equation}
\label{ineq:sigma}
\sigma_0+\sigma_1+\sigma_\infty\geq4.
\end{equation}

We begin with the proof that $m_0=2$.
Suppose that all numbers $m_k\geq3$, then from the master inequality
we deduce,
$$
3n-\sigma_0-\sigma_1-\sigma_\infty\geq 3n-3\quad\Rightarrow\quad
\sigma_0+\sigma_1+\sigma_\infty\leq3,
$$
which contradicts Inequality~(\ref{ineq:sigma}). Now, suppose that $m_1\geq5$.
In this case from the master inequality we get,
\begin{gather*}
\frac{n-\sigma_0}2+\frac{n-\sigma_1}5+\frac{n-\sigma_\infty}5\geq n-1\quad
\Rightarrow\quad10\geq n+5\sigma_0+2\sigma_1+2\sigma_\infty\quad\Rightarrow\\
\quad2\geq n+3\sigma_0\qquad\Rightarrow\qquad\sigma_0=0,\quad n=2.
\end{gather*}
Since $n=2$ the apparent sets in the second and third boxes are empty, thus
the corresponding $R_4(3)$-type reads $R(2|1+1|1+1)$. The latter transformation can be 
treated as belonging to any of the $D$-types mentioned in the Proposition. Explicit form 
of the corresponding $RS_4^2(3)$-trans\-for\-ma\-tion can be found in \cite{AK2} 
(Section~2).

Consider $D$-type $<2,4,m_\infty>$ with $m_\infty\geq5$. the master inequality implies:
\begin{gather*}
\frac{n-\sigma_0}2+\frac{n-\sigma_1}4+\frac{n-\sigma_\infty}5\geq n-1\quad\Rightarrow\quad
20\geq n+10\sigma_0+5\sigma_1+4\sigma_\infty\quad\Rightarrow\\
4\geq n+6\sigma_0+\sigma_1\;\;\Rightarrow\;\;
n=4,\;\sigma_0=\sigma_1=0,\;\sigma_\infty=4,\quad\text{or}\quad
n=2,\;\sigma_0=0,\;\sigma_1=\sigma_\infty=2.
\end{gather*}
The logical case $n=3$ is excluded because it contradicts to the condition $\sigma_0=0$ 
which holds for all $n$. Thus we get two $R_4(3)$-types: $R(2|1+1|1+1)$ and 
$R(2+2|4|\underbrace{1+\ldots+1}_4)$. The last $R$-type has an empty apparent set in the 
last box and, thus can also be treated as belonging to the $D$-type $<2|4|4>$. The rational 
function with this $R$-type exists and the corresponding $RS_4^2(3)$-trans\-for\-ma\-tion 
is explicitly constructed in \cite{AK2} (Section~4, Subsection~4.1.4).

Consider finally $D$-types $<2,3,m_\infty>$ with $m_\infty\geq9$. From the master inequality we find:
\begin{gather}
\frac{n-\sigma_0}2+\frac{n-\sigma_1}3+\frac{n-\sigma_\infty}9\geq n-1\quad\Rightarrow\quad
18\geq n+9\sigma_0+6\sigma_1+2\sigma_\infty\quad\Rightarrow\nonumber\\
10\geq n+7\sigma_0+4\sigma_1\;\;\Rightarrow\;\;
\sigma_0=1,\;\;\sigma_1=0,\;\;n=3,\;\;\sigma_\infty=3,
 \label{eq:sigma10}\\
\sigma_0=0,\quad\sigma_1=0,\quad\;n=2,\ldots,10;
 \label{eq:sigma00}\\
\sigma_0=0,\quad\sigma_1=1,\quad\;n=2,\ldots,6;
 \label{eq:sigma01}\\
\sigma_0=0,\quad\sigma_1=2,\quad\;n=2.
 \label{eq:sigma02}
\end{gather}
In the solution given by Equation~(\ref{eq:sigma10}) we excluded
the case $n=2$, which agrees with the master inequality, because it contradicts to the 
condition $\sigma_0=1$. There is only one $R_4(3)$-type corresponding to 
solution~(\ref{eq:sigma10}), namely $R(2+1|3|1+1+1)$. Because the apparent set in the last 
box is empty this $R$-type can be associated with any $D$-type of the form $<2,3,m>$ with 
arbitrary $m\geq3$, in particular, with $m<9$. The corresponding
rational mapping exists and explicit form of the $RS_4^2(3)$-transformations is given 
in \cite{AK2} (Section~3, Subsection~3.1.2).\\
In the solution~(\ref{eq:sigma00}) $n$ should be divisible by $2$ and $3$, thus $n=6$, 
the apparent set in the last box is empty and hence $\sigma_\infty=6$. There are two 
corresponding $R_4(3)$-types: $R(2+2+2|3+3|2+2+1+1)$ and $R(2+2+2|3+3|3+1+1+1)$. In both 
cases the corresponding rational functions exist, see their explicit forms and corresponding 
solutions of Equation~(\ref{eq:P6}) in \cite{HGBAA} (Section~3, Subsection~3.3, Examples 1 
and 2). Again by the analogous reasoning as in the previous case to both rational functions 
we can assign the same $D$-type $<2,3,m>$, with $m\leq8$. Note that in this case we can
also assign to the first function $D$-type $<2,2,3>$, because in this case we can choose 
the apparent set in the first box consisting of one number $2$, and the nonapparent one - 
with two numbers, both equal $2$, dividing the rest boxes into the apparent and nonapparent 
sets into the natural way we still get the function of $R_4(3)$-type.\\
Turning to the solution~(\ref{eq:sigma01}). We see that $n$ should be even and has the form 
$1+3k$ with some integer $k$. Thus the only possibility is $n=\sigma_\infty=4$ and the 
apparent set in the last box is empty. The only $R_4(3)$-type is $R(2+2|3+1|2+1+1)$. 
The corresponding rational function exists and related $RS_4^2(3)$-trans\-for\-ma\-tion are 
constructed in \cite{AK2} (Section~4, Subsection~4.1.7).\\
Finally, the only $R_4(3)$-type corresponding to Equation~(\ref{eq:sigma02}) is 
$R(2|1+1|1+1)$ is already discussed above.
\end{proof}
\begin{remark}{\rm
To each of the $D$-types, except $<2,3,7>$ and $<2,3,8>$, in Theorem~\ref{thm:main} 
correspond infinite series of rational functions of $R_4(3)$-types. There is a finite 
number of rational functions of $R_4(3)$-type corresponding to the two exceptional 
$D$-types. The latter $D$-types are studied in the subsequent Sections~\ref{sec:D7} and 
\ref{sec:D8}, respectively. It is also not too complicated to describe explicitly the 
infinite series, we plan to do it in further publications.}
\end{remark}
\section{Classification of $RS$-Transformations of $D$-Type $<2,3,7>$}
 \label{sec:D7}
\begin{proposition}
 \label{prop:R237}
There are only three $R_4(3)$-types, with the nonempty apparent set in the third box, 
corresponding to the $D$-type $<2,3,7>$,
namely\footnote{In the numbering of boxes we follow the convention of 
Remark~\ref{rem:convention}.},
\begin{align}
 \label{eq:R7111}
\deg R_4&=10:\quad R_4(7+1+1+1|\underbrace{2+\ldots+2}_5|3+3+3+1),\\
 \label{eq:R72111}
\deg R_4&=12:\quad R_4(7+2+1+1+1|\underbrace{2+\ldots+2}_6
|\underbrace{3+\ldots+3}_4),\\
 \label{eq:R771111}
\deg R_4&=18:\quad R_4(7+7+\underbrace{1+\ldots+1}_4|
\underbrace{2+\ldots+2}_9|\underbrace{3+\ldots+3}_6).
\end{align}
\end{proposition}
\begin{proof}
Put in the master inequality $m_0=2$, $m_1=3$, and $m_\infty=7$, then we
can rewrite it as follows:
$$
42-21\sigma_0-14\sigma_1-6\sigma_\infty\geq n.
$$
Taking into account that $n\geq m_\infty\geq7$ and
Inequality~(\ref{ineq:sigma}) we obtain,
\begin{equation}
 \label{ineq:sigma7}
18-8\sigma_1-15\sigma_0\geq n\geq7.
\end{equation}
The solution of Diophantine Inequality~(\ref{ineq:sigma7}) reads:
\begin{align}
 \label{eq:sigma71}
\sigma_0=0,\quad\sigma_1=0,\quad\sigma_\infty\geq4,\quad 7\leq n\leq18,\\
\label{eq:sigma72}
\sigma_0=0,\quad\sigma_1=1,\quad\sigma_\infty\geq3,\quad 7\leq n\leq10.
\end{align}
Note that Solutions~(\ref{eq:sigma71}) and (\ref{eq:sigma72}) completely
define the second and third boxes of the possible $R$-types.

Consider Solution~(\ref{eq:sigma71}). In this case, $n$ is divisible by
$2\cdot3=6$. Thus the only possibilities are $n=12$ or $n=18$.

If $n=12$ the only possibility is that the apparent set of the first box of
the $R$-type contains only one element. Thus there is only one $R$-type in
this case, namely, (\ref{eq:R72111}). The corresponding covering  and
algebraic solution was constructed in my work \cite{HGBAA} Section 3,
Subsection 3.4, Example 3 (Cross). The same algebraic solution was also
constructed in about the same time by P. Boalch \cite{Bo} by an elaboration
of the method of B.~Dubrovin and M.~Mazzocco \cite{DM}.

If $n=18$ there are two main possibilities:
\begin{enumerate}
\item
The apparent set of the first box consists of two elements the only possible
$R$-type is (\ref{eq:R771111}), because the second and third boxes
are completely defined. Below we show the deformation dessin for
this $R$-type confirming that the corresponding covering really
exists.
\item
The apparent set of the first box consists of one element. There are several logical
possibilities corresponding to the partitions of $18-7=11$ into four natural
numbers. No one of these $R$-types corresponds to a rational covering.
Actually, recall that Euler characteristics of the sphere is $2$,
\begin{equation}
 \label{eq:EULER}
V-E+F=2.
\end{equation}
Suppose that there exists a deformation dessin on the sphere corresponding to 
some of these $R$-types: $V$ is the number of black points plus the
blue one; $F$ is the number of faces which is counted as the four faces,
corresponding to the non-apparent set, plus one face from the apparent set;
and, finally, the valencies of the black points are 3 and valency of the blue one
is 4, each edge is incidental to two vertices:
$$
V=6+1=7,\qquad F=4+1=5.\quad\text{and}\quad E=\frac{3\cdot6+4}2=11.
$$
Now we have $7-11+5=1$, this contradicts Equation~(\ref{eq:EULER}).
\end{enumerate}

Consider now Solution~(\ref{eq:sigma72}). Since $\sigma_0=0$ we
have that $(n|2)>1$, therefore the only logical possibilities are $n=8$,
and $n=10$. For $n=8$ we must have $8=3\cdot k+1$ for some integer $k$,
which is a contradiction. In the case $n=10$ we have
$\sigma_\infty=10-7=3$; together with the facts that $\sigma_0=0$,
$\sigma_1=1$, and that the total number of points in the non-apparent set is $4$,
this implies that there is only one $R$-type corresponding to this case,
namely, (\ref{eq:R7111}). 
\end{proof}

Now we turn to the discussion of existence and explicit constructions of
rational functions with the $R$-types~(\ref{eq:R7111}) and (\ref{eq:R771111}),
as is mentioned in the proof the function with $R$-type~(\ref{eq:R72111}) is already
constructed in \cite{HGBAA}.
 
Consider $R$-type (\ref{eq:R771111}), to confirm the existence of the
corresponding covering we have yet  to present the corresponding
deformation dessin. Note that the type is reducible,
\begin{align}
 \label{eq:R771111-composition}
&R(7+7+\underbrace{1+\ldots+1}_4|\underbrace{2+\ldots+2}_9
|\underbrace{3+\ldots+3}_6)=\\
 \label{eq:R771111-split}
&R(7+1+1|\underbrace{2+\ldots+2}_4+\underset{\wedge}1|3+3+3)\circ
R(1+1|\underset{\wedge}2|1+1)
\end{align}
\begin{remark}
 \label{rem:hyper711}{\rm
This is a digression to the theory of the Gauss hypergeometric functions.
The irreducible Belyi function
$R(7+1+1|\underbrace{2+\ldots+2}_4+1|3+3+3)$ is
of $R_3(3)$-type and defines the following ``seed'' $RS$-transformations
\footnote{The extended notation for $RS$-trans\-for\-ma\-tions that we use below
is explained in \cite{HGBAA,AK1}},
$$
RS_3^2\left(\!\!
\begin{tabular}{c|c|c}
$k/7$&1/2&1/3\\
7+1+1&$\underbrace{2+\!\ldots\!+2}_4+1$&$3+3+3$
\end{tabular}\!\!
\right),\qquad k=1,2,3,
$$
which are equivalent to three seed transformations of the Gauss hypergeometric
functions of order $9$ and in terms of the
$\theta$-triples\footnote{\label{ftn:theta}
$\theta$-triples, the set of formal monodromies for the matrix form of
the hypergeometric equation (see, \cite{AK1,HGBAA}), which differ from
the standard triples of the local exponents for the canonical (scalar)
form of the Gauss hypergeometric equation by the shift $1$ in one of
the elements.} read:
$$
\left(\frac{k}7,\frac12,\frac13\right)\rightarrow
\left(\frac{k}7,\frac{k}7,\frac12\right),\qquad k=1,2,3.
$$
Each of these three transformations allows one to enlarge (by two lines) 
one of the corresponding Octic Clusters introduced in Section~5 of \cite{HGBAA}, because
to the resulted hypergeometric function one can apply a proper quadratic transformation.
}\end{remark}

One of the simplest forms of the first $R$-function in
Equation~(\ref{eq:R771111-split}) is
$$
\lambda_1=\frac{3^{10}(16\lambda^2+7\lambda+49)}
{8(32\lambda^3+84\lambda-35)^3},\;\;
\lambda_1-1=-\frac{(512\lambda^4+256\lambda^3+2208\lambda^2+
328\lambda+1799)^2(\lambda-1)}{8(32\lambda^3+84\lambda-35)^3}.
$$
For application to the theory of algebraic solutions as well as for the
Gauss hypergeometric functions we need the following cumbersome looking
normalizations of this function:
\begin{align}
 \label{eq:norm711-1}
&\lambda_1=
\frac{27(13+7i\sqrt{7})(\lambda-1)
(\lambda-\frac{7^3}{2^8}-\frac{7\cdot13}{2^8}i\sqrt{7})^7}
{2\left(\lambda^3+\frac{7}{2^4}(49+29i\sqrt{7})\lambda^2-
\frac{7^3}{2^{10}}(129+29i\sqrt{7})\lambda+
\frac{3^2\cdot7^3}{2^{15}}(\frac{7^2\cdot13}3-29i\sqrt{7})\right)^3},\\
 \label{eq:norm711-2}
&\lambda_1=\frac{\lambda(\lambda-1)(\lambda-\frac12+\frac{13}{98}i\sqrt{7})^7}
{\left(\lambda^3-(\frac32+\frac{29}{42}i\sqrt{7})\lambda^2+(\frac{121}{294}+
\frac{29}{42}i\sqrt{7})\lambda+\frac{13}{294}+
\frac{29}{4802}i\sqrt{7}\right)^3}
\end{align}
We factorized integers in Equation~(\ref{eq:norm711-1}) only for the purpose of
fitting on one line.
\begin{remark}
 \label{rem:hyper-V}{\rm
While this paper was under preparation I got an information from
R.~Vi\-d${\rm\bar{u}}$\-nas about his recent paper~\cite{V} on
classification of pull-back transformations for the Gauss hypergeometric
functions. This paper is giving a nice and quite profound account of this
subject, in particular, one finds there Equation~(\ref{eq:norm711-2}) in
a slightly different notation. Some of the other Belyi functions of 
$R_3(3)$-type that we discuss in this work were constructed by 
R.~Vi\-d${\rm\bar{u}}$\-nas with the help of the method developed in his 
earlier work~\cite{V1}. The previous Remark~\ref{rem:hyper711} gives also
an illustration to the statement made in Introduction that 
$RS$-transformations seems to be a more general ones than the algebraic
pull-back transformations considered in \cite{V}. Because the local exponent
differences, in the language of \cite{V}, for the $RS$-transformations should
not necessarily be equal to inverse integers as is assumed in \cite{V}:
with each rational covering, in general, we associate a few
{\em independent} (seed) transformations of the Gauss hypergeometric
function, see Remark~\ref{rem:hyper711} above and exact examples in
\cite{AK1}. The situation in this respect is similar with the construction of
algebraic solutions for Equation~(\ref{eq:P6}). There are also some 
intersections of \cite{V} with Sections~4 and~5 work \cite{HGBAA}.}
\end{remark}

To get an explicit realization of Equation~(\ref{eq:R771111-composition}),
~(\ref{eq:R771111-split}) we have to present a rational function
of the $R$-type $R(1+1|21+1)$ in a suitable normalization:
\begin{equation}
 \label{eq:typeR2}
\lambda=\frac{(1-2s)(\lambda_2-t/s)^2}{(\lambda_2-t)},\qquad
t=\frac{s^2}{(2s-1)}.
\end{equation}
Note that the rational function
$\hat\lambda=\hat\lambda(\lambda_2)\equiv\lambda-1$, where
$\lambda$ is given by Equation~(\ref{eq:typeR2}),
is correctly normalized in the sense of Theorem 2.1 \cite{HGBAA}:
$\hat\lambda=\frac{(1-2s)\lambda_2(\lambda_2-1)}{(\lambda_2-t)}$.
Applying this Theorem we calculate algebraic solution of
Equation~(\ref{eq:P6}),
\begin{equation}
 \label{eq:sol2}
t=\frac{s^2}{2s-1}.\qquad y(t)=s=t+\sqrt{t^2-t},
\end{equation}
corresponding to the following $\theta$-tuple,
\begin{equation}
 \label{eq:theta-sol2}
\theta_0=\theta_1,\qquad\theta_t=1-\theta_\infty,
\end{equation}
with two parameters $\theta_0$ and $\theta_t\in\mathbb{C}$.
Substituting $\lambda$ given by Equation~(\ref{eq:typeR2})
into Equation~(\ref{eq:norm711-1}) one gets a rational function of the
$R$-type (\ref{eq:R771111-composition}) correctly normalized in the sense
of Theorem 2.1 of \cite{HGBAA}. Clearly the only critical points of such
composed function which depends on $s$ should coincide with the critical points
of the function (\ref{eq:typeR2}), thus the algebraic solution defined by the
composition exactly coincide with the solution~(\ref{eq:sol2}), however,
now Theorem 2.1 gives for this solution a more restricted $\theta$-tuple:
$\theta_0=\theta_1=\theta_t=1-\theta_\infty=1/7$.

Are there any other rational functions of the
$R$-type~(\ref{eq:R771111})? To answer the question let's
study the following problem: how one can get the functions of this type via
deformations of dessins d'enfants?

It will be convenient to define a notion of {\it symmetric dessins}.
We call a dessin d'enfant or deformation dessin symmetric if it is homeomorphic to
a graph on the Riemann sphere which is invariant under the involution
$\lambda\to-\lambda$. In this case the rational function corresponding to such dessin
can be presented as a composition of a quadratic rational function with a rational
function with the twice lower degree than the original one. In the case of dessins
d'enfants both rational functions are, of course,  the Belyi functions, while for
the deformation dessins the first function of the composition is the Belyi function,
while the second one is a one dimensional deformation of the quadratic Belyi function.
The latter is unique modulo fractional linear transformation of the critical points.
Suppose we consider deformations of a symmetric dessin d'enfant.
If the deformation dessin is symmetric, then we call it the
{\it symmetry preserving deformation}; if the deformation dessin is not symmetric -
the {\it symmetry braking deformations}.
\begin{proposition}
 \label{prop:R771111-symmetric}
Deformation dessins corresponding to the $R$--type~{\rm(\ref{eq:R771111})} can be
obtained only as the symmetry preserving deformations of dessins d'enfants.
\end{proposition}
\begin{proof}
We begin with the ``face'' deformations. There are two types of such
deformations: the cross and join, since the twist can be regarded
as a special case of the join). First consider the cross. Such
deformation is dividing one face of a dessin on two faces and can
increase the black order of a face neighboring with the divided
one (if the latter exists). All in all a dessin before the
cross-deformation should have $5$ faces. We call ``heads'' the
faces with the black order $1$. In case the dessin contains
already four heads its $R$-type can be only
$R(14+\underbrace{1+\ldots+1}_4|\underbrace{2+\ldots+2}_9|
\underbrace{3+\ldots+3}_6)$ (see Figure~\ref{fig:771111-face-2}).
Each head contains on the boundary exactly one black point of
valency $3$ and therefore looks like a balloon on a rope or, as
we say, the head on the ``neck''. In this case the only
possibility is to cross with a chosen neck one of the heads, the necks,
or the edge connecting the ``dumbbells''. The edge and the neck
cannot cross themselves because in this case we have an
``illegal'' deformation which contains a face surrounded with this edge
and, therefore, having its black order equal $0$. Because the dessins 
are located on the
sphere we have only two different possibilities of crossing each
neck or the edge, all in all $(8=2\cdot4)$ variants. One checks
that none of them leads to the right face
distribution~(\ref{eq:R771111-composition}). If the dessin before
the deformation contains exactly three heads, then at least two
of them located in one large face, because we cannot have more
then $5$ faces to get after the deformation $6$ ones. The black
order of face with two heads is at least $9$, so that the
remaining face has the black order $\leq6$. So, after the
deformation the black order of the latter face should be
increased to $7$. The only way of such increase is when one of
the heads entering into it, so that the neck of this head crosses
the boundary of the face. This deformation increases the black
order of the face by $2$. This means that the face distribution
of the dessin before the deformation is $18=10+5+1+1+1$. Such
dessin really exists, but its cross deformation of the type we
discuss leads to the face distribution $18=7+6+2+1+1+1$, see
Figure~\ref{fig:771111-proof}.
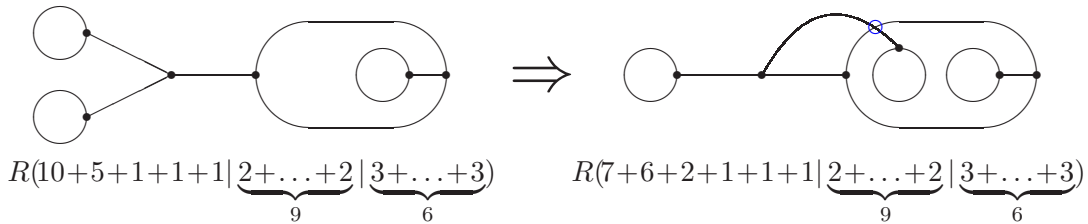
\begin{figure}[ht]
\centering
\begin{picture}(176,72)
\put(20,62){\circle{20}}
\put(130,46){\oval(72,40)}
\put(20,30){\circle{20}}
\put(142,46){\circle{20}}
\put(62,46){\line(1,0){32}}
\put(152,46){\line(1,0){14}}
\put(30,62){\line(2,-1){32}}
\put(30,30){\line(2,+1){32}}
\put(62,46){\circle*{3}}
\put(152,46){\circle*{3}}
\put(166,46){\circle*{3}}
\put(94,46){\circle*{3}}
\put(30,62){\circle*{3}}
\put(30,30){\circle*{3}}
\put(0,06){$R(\!10\!+\!5\!+\!1\!+\!1\!+\!1|\underbrace{2\!+\!\ldots\!+\!2}_9
|\underbrace{3\!+\!\ldots\!+\!3}_6)$}
\end{picture}
\begin{picture}(30,72)
\put(10,40){\Huge$\Rightarrow$}
\end{picture}
\begin{picture}(182,72)
\put(30,46){\circle{20}}
\put(140,46){\oval(72,40)}
\put(124,46){\circle{20}}
\put(152,46){\circle{20}}
\put(72,46){\line(1,0){32}}
\put(162,46){\line(1,0){14}}
\put(40,46){\line(1,0){32}}
\qbezier(72,46)(94,86)(124,56)
\put(72,46){\circle*{3}}
\put(162,46){\circle*{3}}
\put(176,46){\circle*{3}}
\put(104,46){\circle*{3}}
\put(40,46){\circle*{3}}
\put(124,56){\circle*{3}}
\put(115,64){\color{blue}\circle{5}}
\put(0,06){$R(\!7\!+\!6\!+\!2\!+\!1\!+\!1\!+\!1|
\underbrace{2\!+\!\ldots\!+\!2}_9|\underbrace{3\!+\!\ldots\!+\!3}_6)$}
\end{picture}
\caption{An illustration to the proof of non-existence of
cross-deformations of the dessins of $R$-type~(\ref{eq:R771111-composition}).}
\label{fig:771111-proof}
\end{figure}

Suppose now that a dessin before the cross-deformation contains only two
heads. In this case two more heads should appear as a result of the
deformation. The only way it can happen is if the dessin consists of
the circle with one black point on it. The rest of the dessin should
be located inside of the circle and ``live'' on the ``trunk'' which
``grows'' on this black point. If there would be a part of the dessin
outside the circle, then the circle should contain one more black point,
because the valencies of all black points equal $3$.
The deformation in this case is a crossing of the circle by
the trunk such that inside the circle remains only a part of the trunk
while all other parts of the dessin move outside the circle. Because our
pictures are drawn on the Riemann sphere, the dessin before the deformation
actually should contain $3$ heads rather than $2$! In fact, the face
outside the circle where the whole dessin is located is the head.
If the rest of the dessin contains only one more head then, we would have
only three heads as the result of the deformation.

Clearly a dessin before the cross should contain at least two heads, because
there are no one-dimensional deformations that can reduce black orders of
three faces.

Deformation of the join type affects only one face and does not
change the black order of other faces. The affected face is
divided by two ones. So the only possible face distributions of
the dessins that can be deformed by a join to
$R$-type~(\ref{eq:R771111-composition}) are: $18=7+7+2+1+1$,
$18=14+1+1+1+1$, $18=8+7+1+1+1$. The dessin with the last face
distribution does not exist. In fact, suppose that the last dessin 
exists, then it contains three heads. There are two large faces with the black
order equal $8$ and $7$, therefore two heads  are located inside
one of them. Their ``necks'' are connected either with each other
at some point and then this point connected to the boundary of
the surrounding face or with the boundary of the face.
Calculating the black order of such ``minimal construction'' we
get $9$ in the first case and $8$ in the second. However, there
is one more head. This head should be located in the other large
face because, otherwise it cannot have a black order more than
$3$. The last head should be connected with its ``neck'' to the
joint boundary of the large faces at some point different from
the connection points of the other heads, because the valencies
of all connection points equal $3$. Therefore, the minimal black
order of the face containing two heads would be $9$.\\
Figures~\ref{fig:771111-face-1} and ~\ref{fig:771111-face-2} proves that the
first two face deformations really exist.
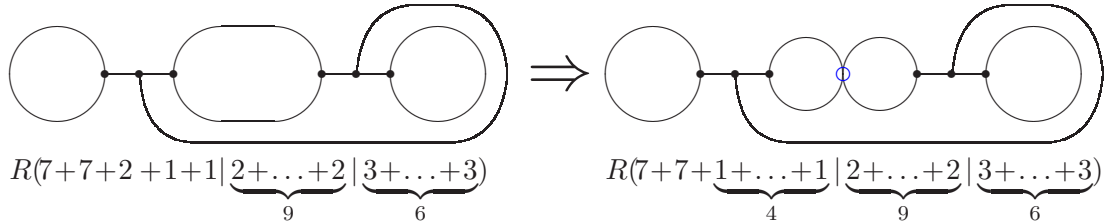
\begin{figure}[ht]
\centering
\begin{picture}(188,74)
\put(18,46){\circle{36}}
\put(90,46){\oval(56,36)}
\put(162,46){\circle{36}}
\put(36,46){\line(1,0){26}}
\put(118,46){\line(1,0){26}}
\put(36,46){\circle*{3}}
\put(49,46){\circle*{3}}
\put(62,46){\circle*{3}}
\put(118,46){\circle*{3}}
\put(131,46){\circle*{3}}
\put(144,46){\circle*{3}}
\qbezier(49,46)(50,21)(69,20)
\put(69,20){\line(1,0){99}}
\qbezier(168,20)(186,22)(188,46)
\qbezier(188,46)(186,70)(168,72)
\qbezier(131,46)(133,70)(151,72)
\put(151,72){\line(1,0){17}}
\put(0,06){$R(\!7\!+\!7\!+\!2+\!1\!+\!1|\underbrace{2\!+\!\ldots\!+\!2}_9
|\underbrace{3\!+\!\ldots\!+\!3}_6)$}
\end{picture}
\begin{picture}(30,74)
\put(04,40){\Huge$\Rightarrow$}
\end{picture}
\begin{picture}(188,74)
\put(18,46){\circle{36}}
\put(162,46){\circle{36}}
\put(36,46){\line(1,0){26}}
\put(118,46){\line(1,0){26}}
\put(36,46){\circle*{3}}
\put(49,46){\circle*{3}}
\put(62,46){\circle*{3}}
\put(118,46){\circle*{3}}
\put(131,46){\circle*{3}}
\put(144,46){\circle*{3}}
\qbezier(49,46)(50,21)(69,20)
\put(69,20){\line(1,0){99}}
\qbezier(168,20)(186,22)(188,46)
\qbezier(188,46)(186,70)(168,72)
\qbezier(131,46)(133,70)(151,72)
\put(151,72){\line(1,0){17}}
\put(104,46){\circle{28}}
\put(76,46){\circle{28}}

\put(90,46){\color{blue}\circle{5}}
\
\put(0,06){$R(\!7\!+\!7\!+\!\underbrace{1\!+\!\ldots\!+\!1}_4
|\underbrace{2\!+\!\ldots\!+\!2}_9|\underbrace{3\!+\!\ldots\!+\!3}_6)$}
\end{picture}
\caption{A symmetry preserving twist of the reducible symmetric dessin.}
\label{fig:771111-face-1}
\end{figure}
\begin{figure}[ht]
\centering
\begin{picture}(176,126)
\put(54,40){\circle{36}}
\put(54,112){\circle{36}}
\put(128,40){\circle{36}}
\put(128,112){\circle{36}}
\put(54,58){\line(0,1){36}}
\put(128,58){\line(0,1){36}}
\put(54,76){\line(1,0){72}}
\put(54,58){\circle*{3}}
\put(128,58){\circle*{3}}
\put(54,76){\circle*{3}}
\put(128,94){\circle*{3}}
\put(54,94){\circle*{3}}
\put(128,76){\circle*{3}}
\put(0,06){$R(\!14\!+\!\underbrace{1\!+\!\ldots\!+\!1}_4|
\underbrace{2\!+\!\ldots\!+\!2}_9|\underbrace{3\!+\!\ldots\!+\!3}_6)$}
\end{picture}
\begin{picture}(30,74)
\put(04,66){\Huge$\Rightarrow$}
\end{picture}
\begin{picture}(188,126)
\put(24,40){\circle{36}}
\put(68,90){\circle{36}}
\put(140,90){\circle{36}}
\put(96,112){\circle{36}}
\qbezier(42,40)(55,43)(68,46)
\put(68,46){\line(0,1){26}}
\put(68,46){\line(1,0){72}}
\put(140,46){\line(0,1){26}}
\put(140,46){\line(1,0){16}}
\qbezier(156,46)(168,48)(170,66)
\put(170,66){\line(0,1){36}}
\qbezier(170,102)(168,118)(156,120)
\put(156,120){\line(-1,0){44}}
\put(42,40){\circle*{3}}
\put(68,72){\circle*{3}}
\put(68,46){\circle*{3}}
\put(112,120){\circle*{3}}
\put(140,46){\circle*{3}}
\put(140,72){\circle*{3}}
\put(82,101){\color{blue}\circle{5}}
\put(0,06){$R(\!7\!+\!7\!+\!\underbrace{1\!+\!\ldots\!+\!1}_4
|\underbrace{2\!+\!\ldots\!+\!2}_9|\underbrace{3\!+\!\ldots\!+\!3}_6)$}
\end{picture}
\caption{A symmetry preserving join of the reducible symmetric dessin.}
\label{fig:771111-face-2}
\end{figure}
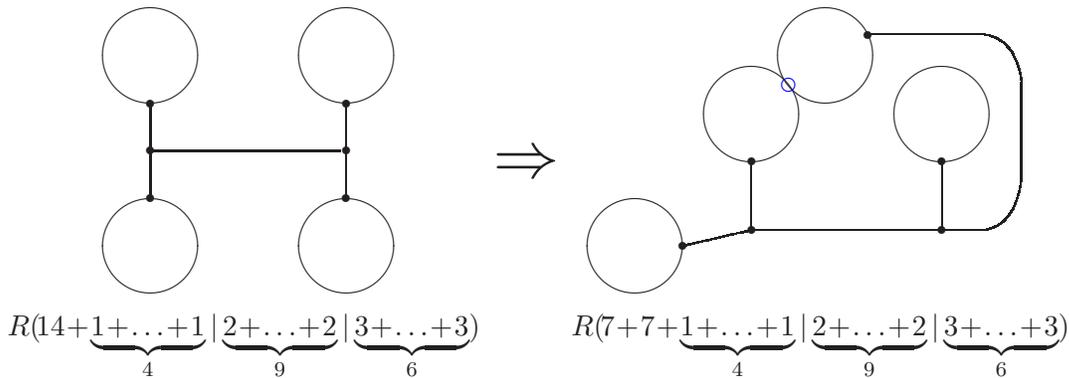

Note that the deformation dessins in the r.-h.s. of
Figures~\ref{fig:771111-face-1} and~\ref{fig:771111-face-2} are
homeomorphic on the Riemann sphere.

Besides the face deformations, there are also deformations which we call ``splits'', 
or, more specifically, $B$- and $W$-splits, depending on the color (black or
white) of the splitted vertex. As we will see below not all such
splits are equivalent. To distinguish different splits we use
notation $LB$- or,say, $CW$-split to denote location of the blue
point after the split, in the first case the blue point belongs
to the crossing of two lines, in the second - of two circles, the last letter
means, of course, the color of the splitted vertex. If the blue
vertex belongs to a circle and line we denote such deformation as
$CLB$-split, if $B$-vertex is splitted.

In our case we obviously have only two splits: $W$-split (4=2+2) and
$B$-split (6=3+3). These deformations are shown on
Figures~\ref{fig:771111-split-W} and~\ref{fig:771111-split-B}.
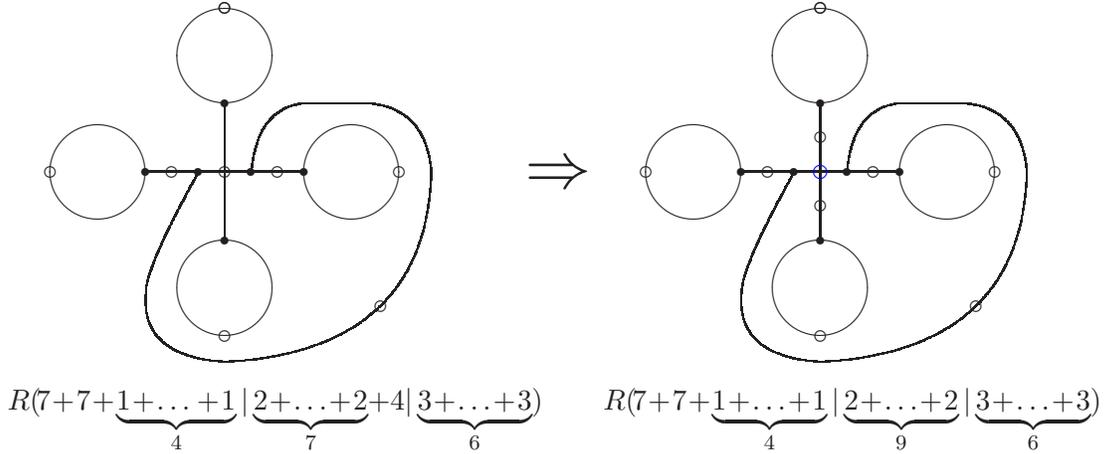
\begin{figure}[ht]
\centering
\begin{picture}(188,158)
\put(34,96){\circle{36}}
\put(130,96){\circle{36}}
\put(82,140){\circle{36}}
\put(82,52){\circle{36}}
\put(52,96){\line(1,0){60}}
\put(82,96){\line(0,1){26}}
\put(82,70){\line(0,1){26}}
\put(82,96){\circle{4}}
\put(52,96){\circle*{3}}
\put(72,96){\circle*{3}}
\put(62,96){\circle{4}}
\put(92,96){\circle*{3}}
\put(102,96){\circle{4}}
\put(112,96){\circle*{3}}
\put(16,96){\circle{4}}
\put(148,96){\circle{4}}
\put(82,158){\circle{4}}
\put(82,122){\circle*{3}}
\put(82,158){\circle{4}}
\put(82,70){\circle*{3}}
\put(82,34){\circle{4}}
\put(141,45){\circle{4}}
\qbezier(72,96)(52,60)(52,50)
\qbezier(52,50)(50,26)(82,24)
\qbezier(82,24)(158,28)(160,96)
\qbezier(160,96)(158,120)(136,122)
\qbezier(92,96)(94,120)(112,122)
\put(112,122){\line(1,0){24}}
\put(0,06){$R(\!7\!+\!7\!+\!\underbrace{1\!+\!\ldots+\!1}_4|
\underbrace{2\!+\!\ldots\!+\!2}_7\!+4|\underbrace{3\!+\!\ldots\!+\!3}_6)$}
\end{picture}
\begin{picture}(30,158)
\put(04,90){\Huge$\Rightarrow$}
\end{picture}
\begin{picture}(188,158)
\put(34,96){\circle{36}}
\put(130,96){\circle{36}}
\put(82,140){\circle{36}}
\put(82,52){\circle{36}}
\put(52,96){\line(1,0){60}}
\put(82,96){\line(0,1){26}}
\put(82,70){\line(0,1){26}}
\put(82,109){\circle{4}}
\put(82,96){\color{blue}\circle{5}}
\put(82,83){\circle{4}}
\put(52,96){\circle*{3}}
\put(72,96){\circle*{3}}
\put(62,96){\circle{4}}
\put(92,96){\circle*{3}}
\put(102,96){\circle{4}}
\put(112,96){\circle*{3}}
\put(16,96){\circle{4}}
\put(148,96){\circle{4}}
\put(82,158){\circle{4}}
\put(82,122){\circle*{3}}
\put(82,158){\circle{4}}
\put(82,70){\circle*{3}}
\put(82,34){\circle{4}}
\put(141,45){\circle{4}}
\qbezier(72,96)(52,60)(52,50)
\qbezier(52,50)(50,26)(82,24)
\qbezier(82,24)(158,28)(160,96)
\qbezier(160,96)(158,120)(136,122)
\qbezier(92,96)(94,120)(112,122)
\put(112,122){\line(1,0){24}}
\put(0,06){$R(\!7\!+\!7\!+\!\underbrace{1\!+\!\ldots\!+\!1}_4
|\underbrace{2\!+\!\ldots\!+\!2}_9|\underbrace{3\!+\!\ldots\!+\!3}_6)$}
\end{picture}
\caption{A symmetry preserving $W$-split of the reducible symmetric dessin.}
\label{fig:771111-split-W}
\end{figure}
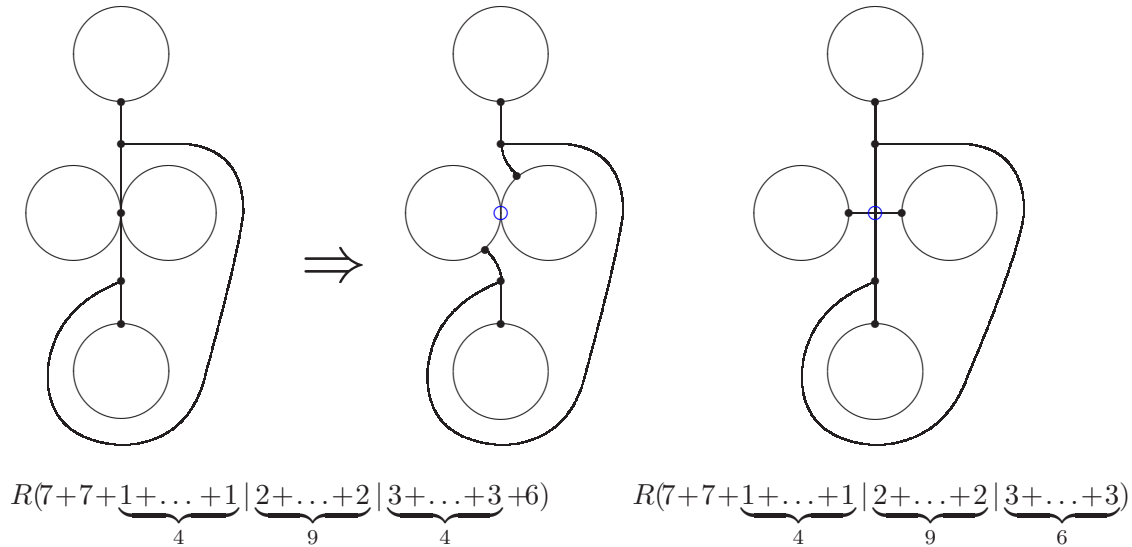
\begin{figure}[ht]
\centering
\begin{picture}(140,194)
\put(24,116){\circle{36}}
\put(60,116){\circle{36}}
\put(42,176){\circle{36}}
\put(42,56){\circle{36}}
\put(42,116){\line(0,1){42}}
\put(42,74){\line(0,1){42}}
\put(42,158){\circle*{3}}
\put(42,116){\circle*{3}}
\put(42,74){\circle*{3}}
\put(42,142){\circle*{3}}
\put(42,90){\circle*{3}}
\qbezier(42,90)(16,80)(14,56)
\qbezier(14,56)(13,30)(42,28)
\qbezier(42,28)(68,29)(74,56)
\qbezier(88,116)(88,140)(66,142)
\qbezier(74,56)(86,100)(88,116)
\put(42,142){\line(1,0){24}}
\put(110,90){\Huge$\Rightarrow$}
\put(0,06){$R(\!7\!+\!7\!+\!\underbrace{1\!+\!\ldots+\!1}_4|
\underbrace{2\!+\!\ldots\!+\!2}_9|\underbrace{3\!+\!\ldots\!+\!3}_4+\!6)$}
\end{picture}
\begin{picture}(88,194)
\put(24,116){\circle{36}}
\put(60,116){\circle{36}}
\put(42,176){\circle{36}}
\put(42,56){\circle{36}}
\put(42,142){\line(0,1){16}}
\qbezier(42,142)(42,136)(48,130)
\put(42,74){\line(0,1){16}}
\qbezier(42,90)(42,96)(36,102)
\put(42,158){\circle*{3}}
\put(36,102){\circle*{3}}
\put(48,130){\circle*{3}}
\put(42,74){\circle*{3}}
\put(42,116){\color{blue}\circle{5}}
\put(42,142){\circle*{3}}
\put(42,90){\circle*{3}}
\qbezier(42,90)(16,80)(14,56)
\qbezier(14,56)(13,30)(42,28)
\qbezier(42,28)(68,29)(74,56)
\qbezier(88,116)(88,140)(66,142)
\qbezier(74,56)(86,100)(88,116)
\put(42,142){\line(1,0){24}}
\end{picture}
\begin{picture}(186,194)
\put(64,116){\circle{36}}
\put(120,116){\circle{36}}
\put(92,176){\circle{36}}
\put(92,56){\circle{36}}
\put(82,116){\line(1,0){20}}
\put(92,116){\line(0,1){42}}
\put(92,74){\line(0,1){42}}
\put(82,116){\circle*{3}}
\put(102,116){\circle*{3}}
\put(92,158){\circle*{3}}
\put(92,116){\color{blue}\circle{5}}
\put(92,74){\circle*{3}}
\put(92,142){\circle*{3}}
\put(92,90){\circle*{3}}
\qbezier(92,90)(66,80)(64,56)
\qbezier(64,56)(63,30)(92,28)
\qbezier(92,28)(118,29)(127,55)
\qbezier(148,116)(148,140)(126,142)
\qbezier(127,55)(146,100)(148,116)
\put(92,142){\line(1,0){34}}
\put(0,06){$R(\!7\!+\!7\!+\!\underbrace{1\!+\!\ldots\!+\!1}_4
|\underbrace{2\!+\!\ldots\!+\!2}_9|\underbrace{3\!+\!\ldots\!+\!3}_6)$}
\end{picture}
\caption{Symmetry preserving $CB$- and $LB$-splits of the reducible
symmetric dessin.}
\label{fig:771111-split-B}
\end{figure}

Note that $W$-split on Figures~\ref{fig:771111-split-W} is homeomorphic
in the Riemann sphere to $LB$-split on Figure~\ref{fig:771111-split-B}.
Also $CB$-split on Figure~\ref{fig:771111-split-B} is
homeomorphic on the Riemann sphere to the twist on
Figures~\ref{fig:771111-face-1}. Moreover, both deformation dessins on
Figure~\ref{fig:771111-split-B} represent two branches of the same
rational covering, because, clearly, they are homotopic, continuously deformable 
one into another through the dessin in the l.-h.s. of this picture.

There is also $CLB$-split of the dessin on
Figure~\ref{fig:771111-split-B} with the right valencies of black
(and, of course, white) vertices (see Figure~\ref{fig:771111-split-CLB}).
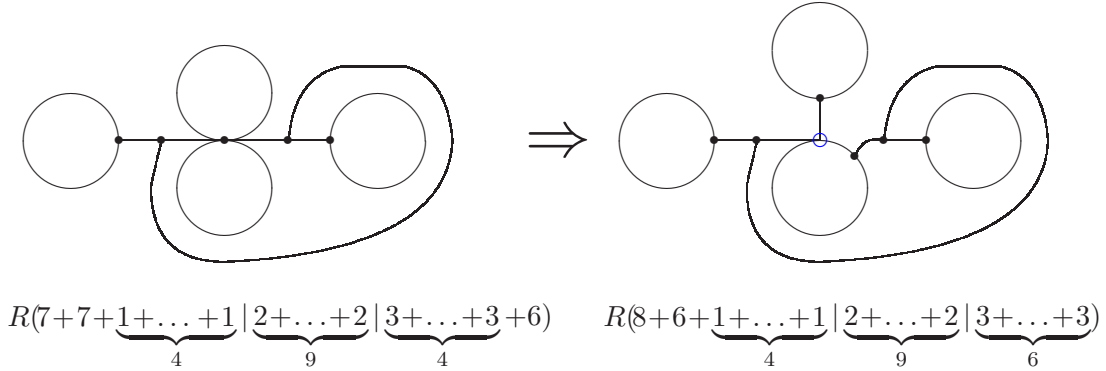
\begin{figure}[ht]
\centering
\begin{picture}(188,128)
\put(24,76){\circle{36}}
\put(140,76){\circle{36}}
\put(82,94){\circle{36}}
\put(82,58){\circle{36}}
\put(42,76){\line(1,0){80}}
\put(42,76){\circle*{3}}
\put(58,76){\circle*{3}}
\put(82,76){\circle*{3}}
\put(106,76){\circle*{3}}
\put(122,76){\circle*{3}}
\qbezier(58,76)(54,60)(54,58)
\qbezier(54,58)(54,30)(82,30)
\qbezier(168,76)(166,34)(82,30)
\qbezier(168,76)(166,100)(150,104)
\qbezier(106,76)(108,100)(126,104)
\put(126,104){\line(1,0){24}}
\put(0,06){$R(\!7\!+\!7\!+\!\underbrace{1\!+\!\ldots+\!1}_4|
\underbrace{2\!+\!\ldots\!+\!2}_9|\underbrace{3\!+\!\ldots\!+\!3}_4+6)$}
\end{picture}
\begin{picture}(30,128)
\put(04,70){\Huge$\Rightarrow$}
\end{picture}
\begin{picture}(188,128)
\put(24,76){\circle{36}}
\put(140,76){\circle{36}}
\put(82,110){\circle{36}}
\put(82,58){\circle{36}}
\put(82,76){\line(0,1){16}}
\put(42,76){\line(1,0){40}}
\put(106,76){\line(1,0){16}}
\qbezier(106,76)(98,78)(95,70)
\put(95,70){\circle*{3}}
\put(42,76){\circle*{3}}
\put(58,76){\circle*{3}}
\put(82,92){\circle*{3}}
\put(82,76){\color{blue}\circle{5}}
\put(106,76){\circle*{3}}
\put(122,76){\circle*{3}}
\qbezier(58,76)(54,60)(54,58)
\qbezier(54,58)(54,30)(82,30)
\qbezier(168,76)(166,34)(82,30)
\qbezier(168,76)(166,100)(150,104)
\qbezier(106,76)(108,100)(126,104)
\put(126,104){\line(1,0){24}}
\put(0,06){$R(\!8\!+\!6\!+\!\underbrace{1\!+\!\ldots\!+\!1}_4
|\underbrace{2\!+\!\ldots\!+\!2}_9|\underbrace{3\!+\!\ldots\!+\!3}_6)$}
\end{picture}
\caption{``Symmetry braking''  $CLB$-split of the reducible symmetric dessin.}
\label{fig:771111-split-CLB}
\end{figure}
However the resulted deformation dessin does not belong to $R_4$-type.
\end{proof}

Finally, consider $R$-type~(\ref{eq:R7111}) of Proposition~\ref{prop:R237}.
It can be obtained as: (1) face deformations of the following dessins,
$R(\!8\!+\!1\!+\!1|\underbrace{2\!+\!\ldots\!+\!2}_5|3\!+\!3\!+\!3\!+\!1)$ and
$R(7\!+\!2\!+\!1|\underbrace{2\!+\!\ldots\!+\!2}_5|3\!+\!3\!+\!3\!+\!1)$;
(2) $W$- split of
$R(7\!+\!1\!+\!1\!+\!1|4+\!\underbrace{2\!+\!\ldots\!+\!2}_4|3\!+\!3\!+\!3\!+\!1)$;
or (3) $B$-splits of $R(7+\!1\!+\!1\!+\!1|\underbrace{2\!+\!\ldots\!+\!2}_5|3\!+\!3\!+\!4)$
and $R(7+\!1\!+\!1\!+\!1|\underbrace{2\!+\!\ldots\!+\!2}_5|6\!+\!3\!+\!1)$.
We leave to the interested reader to prove that all these dessins are homotopic in
the Riemann sphere so that the corresponding deformation dessins represent different
branches of one and the same algebraic function. Instead of studying the dessins we
present below an explicit form of the corresponding rational covering: 
\begin{align}
 \label{eq:7111}
 \lambda_1&=\frac{1728s^{12}(3s^2-4s+4)}{(s+2)^{14}(s-1)^8}
 \frac{(\lambda_2-1)(\lambda_2^2+a_1\lambda_2+a_0)}
 {\lambda_2(\lambda_2^3+c_2\lambda_2^2+c_1\lambda_2+c_0)^3},\qquad
 \lambda_2=\frac{\lambda A}{\lambda-1+A},\\
\intertext{where}
a_0&=\frac{27s^4(2s^2-3s+2)^2}{(s+2)^4(s-1)^4(3s^2-4s+4)},\quad
a_1=-\frac{2(14s^5-25s^4+20s^3+8s^2-16s+8)}{(s+2)^2(s-1)^2(3s^2-4s+4)},\nonumber\\
c_0&=-\frac{24s^4(4s^3-s^2-4s+4)}{(s+2)^6(s-1)^4},\qquad
c_1=\frac{60s^6-84s^5-15s^4+72s^3-8s^2-32s+16}{(s+2)^4(s-1)^4},\nonumber\\
c_2&=-\frac{2(6s^3-3s^2-4s+4)}{(s+2)^2(s-1)^2},\qquad\text{and}\nonumber\\
A&=\frac{14s^5-25s^4+20s^3+8(s-1)^2+8(s-1)(s^2-s+1)\sqrt{(2s+1)(1-s)(s^2-s+1)}}
{(s+2)^2(s-1)^2(3s^2-4s+4)},\nonumber
\end{align}
is a solution of the quadratic equation, $\lambda_2^2+a_1\lambda_2+a_0=0$.
Note that the function $\lambda_1=\lambda_1(\lambda_2)$ has a rational
parametrization, however it is not correctly normalized.
After a normalization, the fractional-linear transformation
$\lambda_2=\lambda_2(\lambda)$, we get the function $\lambda_1=\lambda_1(\lambda)$,
which has an elliptic parametrization. Applying now Theorem~2.1 of \cite{HGBAA}, we
get an algebraic solution of Equation~(\ref{eq:P6}),
\begin{align}
 \label{eq:y7111}
&y(t)=1+\frac{(3s-2)(s^2-2s+4)^2}{4(s+2)(s-1)^2(s^2-s+1)(3s^2-4s+4)}\times\\
&\frac{-14s^5+25s^4-20s^3-8s^2+16s-8-8(s-1)(s^2-s+1)
\sqrt{(2s+1)(1-s)(s^2-s+1)}}
{(2s+1)(3s^3-10s^2+6s-2)-14(s-1)\sqrt{(2s+1)(1-s)(s^2-s+1)}}\nonumber\\
 \label{eq:t7111}
&t=\frac12-\sum
\frac{14s^9-105s^8+252s^7-392s^6+420s^5-336s^4+112s^3+72s^2-96s+32}
{16(s+2)^2(s-1)^3(s^2-s+1)\sqrt{(2s+1)(1-s)(s^2-s+1)}},
\end{align}
for the following set of $\theta$-parameters:
$$
\theta_0=\frac13,\qquad\theta_1=\frac17,\qquad
\theta_t=\frac17,\qquad\theta_\infty=\frac67.
$$

There are a few other suitable normalizations of
the function $\lambda_1(\lambda_2)$, clearly all of them can be parameterized only
by algebraic curves of genus $1$. Theorem~2.1~\cite{HGBAA} allows to find an algebraic
solution (of genus 1) to each such normalization. However, it is easy to check that
all these solutions are related to each other via so-called B\"acklund transformations
for Equation~(\ref{eq:P6}). Thus, Equations~(\ref{eq:y7111}) and (\ref{eq:t7111})
represent the only "pull-back" seed algebraic solution. The list of the "$RS$" seed
algebraic solutions is given below in Proposition~\ref{prop:237-RS-7111}. 
We can summarize our study as the following Propositions.
\begin{proposition}
 \label{prop:237-existence}
For all $R$-types specified in Proposition~{\rm\ref{prop:R237}}: {\rm(\ref{eq:R7111}),
(\ref{eq:R72111})}, and {\rm(\ref{eq:R771111})}, there exist rational functions with
these $R$-types. Each of these rational functions can be rationally parameterized by 
a "deformation" parameter $s\in\mathbb{CP}^1\setminus\mathcal{B}$, where $\mathcal{B}$ 
is a finite set. The resulting birational functions: $\lambda_1=\lambda_1(\lambda,s)$
{\rm(}Equation~{\rm(\ref{eq:7111}))}, $\lambda_1=\lambda_1(\lambda_2,s)$
{\rm(}Equations~{\rm(\ref{eq:typeR2})} and {\rm(\ref{eq:norm711-1})}, and
$z=z(z_1,s)$ in {\rm\cite{HGBAA}} Section {\rm3}, Subsection {\rm3.4},
Example {\rm3} {\rm(}Cross{\rm)}, are unique up to fractional-linear transformations of 
the first argument and reparametrization of $s$.
\end{proposition}
\begin{proposition}
 \label{prop:237-RS-7111}
There are three seed $RS$-trans\-for\-mations related with
$R$-type {\rm(\ref{eq:R7111})}:
$
RS_4^2\left(\!\!
\begin{tabular}{c|c|c}
$k/7$&$1/2$&$1/3$\\
$7+1+1+1$&$\underbrace{2+\!\ldots\!+2}_5$&
$3+3+3+1$
\end{tabular}\!\!
\right)
$
for $k=1,2$, and $3$. Each of these transformations produces one algebraic genus $1$
solution for the following sets of the $\theta$-parameters:
\begin{align*}
k&=1,\qquad&\theta_0&=\frac13,&\theta_1&=\frac17,&\theta_1&=\frac17,
&\theta_\infty&=\frac67,\\
k&=2,\qquad&\theta_0&=\frac13,&\theta_1&=\frac27,&\theta_1&=\frac27,
&\theta_\infty&=\frac27,\\
k&=3,\qquad&\theta_0&=\frac13,&\theta_1&=\frac37,&\theta_1&=\frac37,
&\theta_\infty&=\frac47.
\end{align*}
\end{proposition}
\begin{proposition}
 \label{prop:237-RS-72111}
There are four seed $RS$-trans\-for\-mations related with
$R$-type {\rm(\ref{eq:R72111})}:\\
$
RS_4^2\left(\!\!
\begin{tabular}{c|c|c}
$k/7$&$1/2$&$1/3$\\
$7+2+1+1+1$&$\underbrace{2+\!\ldots\!+2}_6$&
$\underbrace{3+\!\ldots\!+3}_4$
\end{tabular}\!\!
\right)
$
for $k=1,2,3$, and $7/2$.
Each of these transformations produces one algebraic genus $0$
solution for the following sets of the $\theta$-parameters:
\begin{align*}
k&=1,\qquad&\theta_0&=\frac17,&\theta_1&=\frac17,&\theta_1&=\frac17,
&\theta_\infty&=\frac57,\\
k&=2,\qquad&\theta_0&=\frac27,&\theta_1&=\frac27,&\theta_1&=\frac27,
&\theta_\infty&=\frac47,\\
k&=3,\qquad&\theta_0&=\frac37,&\theta_1&=\frac37,&\theta_1&=\frac37,
&\theta_\infty&=\frac17,\\
k&=\frac72,\qquad&\theta_0&=\frac12,&\theta_1&=\frac12,&\theta_1&=\frac12,
&\theta_\infty&=-\frac52.
\end{align*}
\end{proposition}
\section{Classification of $RS$-Transformations of $D$-Type $<2,3,8>$}
 \label{sec:D8}
\begin{proposition}
There is only one $R_4(3)$-type, with the nonempty apparent set in the third box 
\footnote{In the numbering of boxes we follow the convention of 
Remark~\ref{rem:convention}. Here the third box comes first in the
natural counting.},
corresponding to the $D$-type $<2,3,8>$, namely, 
\begin{equation}
 \label{eq:R81111}
R_4(8+\underbrace{1+\ldots+1}_4|\underbrace{2+\ldots+2}_6
|\underbrace{3+\ldots+3}_4)
\end{equation}
\end{proposition}
\begin{proof}
We again refer to Master Inequality~(\ref{ineq:master}). For our
particular divisors after simple manipulations it can be rewritten as follows:
$$
24-8\sigma_0-5\sigma_1-3(\sigma_0+\sigma_1+\sigma_\infty)\geq n.
$$
Now, taking into account Inequality~(\ref{ineq:sigma}) and the
fact that $n\geq8$, we find that $12-8\sigma_0-5\sigma_1\geq8$.
Therefore, $\sigma_0=\sigma_1=0$ and thus $\sigma_\infty\geq4$.
Again returning to Master Inequality~(\ref{ineq:master}) and
substituting in it $\sigma_0=\sigma_1=0$, we obtain,
$24-3\sigma_\infty\geq n$ which implies that $n\leq12$. On the
other hand $n\geq m_\infty+\sigma_\infty\geq12$. Therefore,
$n=12$, the only possible $R_4(3)$-type with four non-apparent entries
and non-empty apparent set in the third box is equivalent to (\ref{eq:R81111}).
\end{proof}
Because the degree of the function~(\ref{eq:R81111}) is
$12=2\cdot6=3\cdot4=4\cdot3=6\cdot2$, we have to examine whether this
rational function can be presented as a composition of rational functions of
the lower degree. Clearly, that one of these functions should be the Belyi
function and the other a one-dimensional deformation of (another) Belyi
function. The latter generates an algebraic solution of the sixth Painlev\'e
equation. It is easy to see that such composition defines exactly
the same algebraic solution of P6 as its member, the deformed Belyi
function.

In the above factorizations of $12$ into the divisors we assume that the
first function is the Belyi one, while the second is a deformation, therefore
in this sense these decompositions are not commutative. By a
straightforward analysis, just an examination of a few possibilities, we find
that there is actually only one such composition (see also
Remark~\ref{rem:81111i}) below) corresponding to the
factorization $12=6\cdot2$, namely,
\begin{equation}
 \label{eq:composition}
R(8+\!\underbrace{1+\!\ldots\!+1}_4|\underbrace{2+\!\ldots\!+2}_{6}
|\underbrace{3+\!\ldots\!+3}_{4})=R(\underset{\wedge}4+1+1|2+2+2|3+3)
\circ R(\underset{\wedge}2|1+1|1+1).
\end{equation}
The function $R(4+1+1|2+2+2|3+3)$ itself is also reducible,
$R(4+1+1|2+2+2|3+3)=R((\underset{\wedge}2+1|2+\underset{\wedge}1|3)\circ
R(\underset{\wedge}2|\underset{\wedge}2|1+1)$, however it is not important
in the following.
Explicit form of the functions in the r.-h.s. of
Equation~(\ref{eq:composition}) is as follows:
\begin{align}
&\lambda_2=\frac{108\lambda_1^4(\lambda_1-1)}
{(\lambda_1^2-16\lambda_1+16)^3},\qquad
\lambda_1= (1-2s)\frac{(\lambda-t/s)^2}{(\lambda-t)},\nonumber
\intertext{where}
 \label{eq:t-y-2}
&t=\frac{s^2}{2s-1},\quad\text{and}\quad y(t)=s\\
\intertext{is the solution of Equation~(\ref{eq:P6}) for the $\theta$-tuple,}
 \label{eq:theta-2}
&\theta_0=\theta_1,\qquad\theta_t=\theta_\infty-1,
\end{align}
for arbitrary $\theta_0$ and $\theta_\infty\in\mathbb{C}$
(see Theorem 2.1 of \cite{HGBAA}). The r.-h.s. now is easy to find,
\begin{gather}
 \label{eq:lambda81111}
\lambda_1=\frac{108\lambda(\lambda-1)(\lambda-t)(\lambda-t/s)^8}
{(2s-1)(\lambda^4+c_3\lambda^3+c_2\lambda^2+c_1\lambda+c_0)^3},\\
c_3=-\frac{4(s-4)}{(2s-1)},\quad
c_2=-\frac{2(5s^2+16s-8)}{(2s-1)^2},\quad
c_1=\frac{4(7s-4)s^2}{(2s-1)^3},\quad
c_0=\frac{s^4}{(2s-1)^4},\nonumber
\end{gather}
where $t$ is the same as in (\ref{eq:t-y-2}). Applying Theorem 2.1 of
\cite{HGBAA} we again arrive at the solution $y(t)$ defined in
Equations~(\ref{eq:t-y-2}) but now for a particular choice of the
$\theta$-tuple~(\ref{eq:theta-2}), $\theta_0=\theta_t=1/8$.

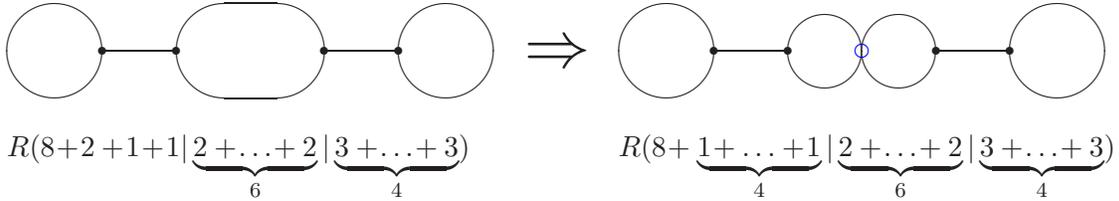
\begin{figure}[ht]
\centering
\begin{picture}(184,70)
\put(18,46){\circle{36}}
\put(92,46){\oval(56,36)}
\put(166,46){\circle{36}}
\put(36,46){\line(1,0){28}}
\put(120,46){\line(1,0){28}}
\put(36,46){\circle*{3}}
\put(64,46){\circle*{3}}
\put(120,46){\circle*{3}}
\put(148,46){\circle*{3}}
\put(0,06){$R(8\!+\!2+\!1\!+\!1|\underbrace{2+\!\ldots\!+2}_6
|\underbrace{3+\!\ldots\!+3}_4)$}
\end{picture}
\begin{picture}(40,70)
\put(08,40){\Huge$\Rightarrow$}
\end{picture}
\begin{picture}(184,70)
\put(18,46){\circle{36}}
\put(106,46){\circle{28}}
\put(78,46){\circle{28}}
\put(166,46){\circle{36}}
\put(36,46){\line(1,0){28}}
\put(120,46){\line(1,0){28}}
\put(92,46){\color{blue}\circle{5}}
\put(36,46){\circle*{3}}
\put(64,46){\circle*{3}}
\put(120,46){\circle*{3}}
\put(148,46){\circle*{3}}
\put(0,06){$R(8\!+\underbrace{1\!+\ldots+\!1}_4|\underbrace{2+\!\ldots\!+2}_6
|\underbrace{3+\!\ldots\!+3}_4)$}
\end{picture}
\caption{A symmetry preserving twist of the reducible dessin for the Belyi
function}
\label{fig:81111-twist}
\end{figure}

It is instructive to confirm the above mentioned analysis that leads to
Equation~(\ref{eq:composition}) graphically
(see Figure~\ref{fig:81111-twist}).

There are two more reducible dessins with the symmetry preserving deformations:
\begin{gather*}
R(8+\!\underbrace{1+\!\ldots\!+1}_4|\underbrace{2+\!\ldots\!+2}_{4}\!+4|
\underbrace{3+\!\ldots\!+3}_{4})=
R(\underset{\wedge}4+1+1|2+2+\underset{\wedge}2|3+3)
\circ R(\underset{\wedge}2|\underset{\wedge}2|1+1),\\
R(8+\underbrace{1+\ldots+1}_4|\underbrace{2+\ldots+2}_{6}|3+3+6)=
R(\underset{\wedge}4+1+1|2+2+2|
3+\underset{\wedge}3)\circ R(\underset{\wedge}2|1+1|\underset{\wedge}2).
\end{gather*}
Their deformation dessins exactly coincide and correspond to another
branch of the solution~(\ref{eq:t-y-2}).

The covering constructed above is not the only possible for this $R$-type.
To get an idea why there should be another solution, we observe that there
is the following deformation:\\
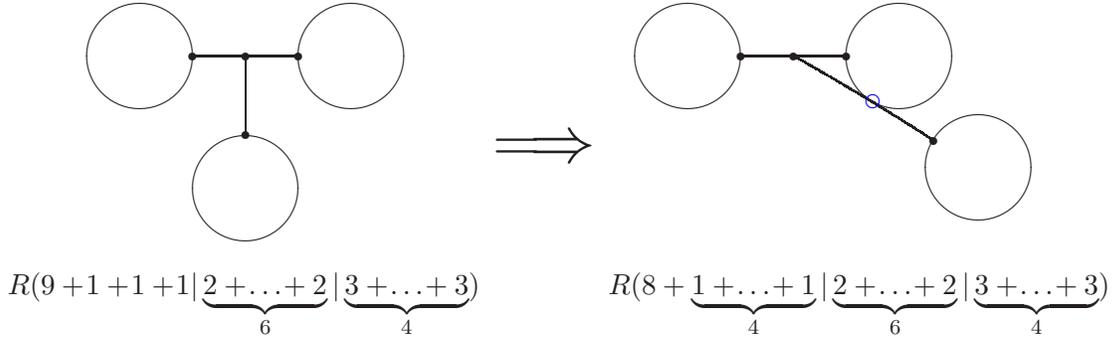
\begin{figure}[ht]
\centering
\begin{picture}(140,120)
\put(30,95){\circle{40}}
\put(110,95){\circle{40}}
\put(50,95){\line(1,0){40}}
\put(50,95){\circle*{3}}
\put(90,95){\circle*{3}}
\put(70,95){\circle*{3}}
\put(70,95){\line(0,-1){30}}
\put(70,45){\circle{40}}
\put(70,65){\circle*{3}}
\put(-20,05){$R(9+\!1+\!1+\!1|\underbrace{2+\!\ldots\!+2}_6|
\underbrace{3+\!\ldots\!+3}_4)$}
\end{picture}
\begin{picture}(60,120)
\put(20,55){\Huge$\Longrightarrow$}
\end{picture}
\begin{picture}(170,90)
\put(30,95){\circle{40}}
\put(110,95){\circle{40}}
\put(50,95){\line(1,0){40}}
\put(50,95){\circle*{3}}
\put(90,95){\circle*{3}}
\put(70,95){\circle*{3}}
\qbezier(70,95)(100,78)(123,63)
\put(100,78){\color{blue}\circle{5}}
\put(140,53){\circle{40}}
\put(123,63){\circle*{3}}
\put(0,05){$R(8+\underbrace{1+\!\ldots\!+1}_4|
\underbrace{2+\!\ldots\!+2}_6|\underbrace{3+\!\ldots\!+3}_4)$}
\end{picture}
\caption{A symmetry braking ``face'' deformation ($LC$-Join) of the dessin for
the reducible Belyi function}
\label{fig:81111i}
\end{figure}
\begin{remark}
 \label{rem:huper9111}{\rm
This is a digression to the theory of the Gauss hypergeometric functions.
The Belyi function in the r.-h.s. of Figure~\ref{fig:81111i},
as well as the Belyi function from Figure~\ref{fig:81111-twist},
have in our terminology $R_3$-type and therefore define transformations of
order $12$ for the Gauss hypergeometric function. Both of
these functions are reducible:
\begin{gather}
 \label{eq:9111-decomposition}
R(9+\!1+\!1+\!1|\underbrace{2+\!\ldots\!+2}_6|
\underbrace{3+\!\ldots\!+3}_4)=
R(\underset{\wedge}3+1|2+2|3+\underset{\wedge}1)\circ
R(\underset{\wedge}3|1+1+1|\underset{\wedge}3),\\
 \label{eq:8211-decomposition}
R(8+2+1+1|\underbrace{2+\!\ldots\!+2}_6
|\underbrace{3+\!\ldots\!+3}_4)=R(2+1|2+1|3)\circ R(2|2|1+1)\circ R(2|1+1|2).
\end{gather}
In most cases reducibility of $R$-function means that the corresponding
higher order transformation for the Gauss hypergeometric function is a
composition of transformations of the lower order, namely those that
correspond to the $R$-functions of lower degree from the decomposition of
the original $R$-function. This happens
when all rational functions in the corresponding decomposition have $R_3$-type.
Like it is in the case of the function~(\ref{eq:8211-decomposition}),
for $k=1,3$:
\begin{gather}
\label{8211rs}
RS_3^2\left(\!\!
\begin{tabular}{c|c|c}
$k/8$&1/2&1/3\\
8+2+1+1&$\underbrace{2+\!\ldots\!+2}_6$&
$\underbrace{3+\!\ldots\!+3}_4$
\end{tabular}\!\!
\right)=\\
RS_3^2\left(\!\!
\begin{tabular}{c|c|c}
$k/8$&1/2&1/3\\
2+1&2+1&3
\end{tabular}\!\!
\right)\circ
RS_3^2\left(\!\!
\begin{tabular}{c|c|c}
$k/4$&1/2&$k/8$\\
2&2&1+1
\end{tabular}\!\!
\right)\circ
RS_3^2\left(\!\!
\begin{tabular}{c|c|c}
$k/2$&$k/8$&$k/8$\\
2&2&1+1
\end{tabular}\!\!
\right)\label{eq:8211rs-decomposition}
\end{gather}
The last transformation is possible because in this case $k/2=1/2\mod(1)$.
The resulting $\theta$-triple\footnote{see Footnote~\ref{ftn:theta} on
page~\pageref{ftn:theta}}
is $\big(k/8,k/8,k/4-1\big)$.
The other two ``seed'' transformations~(\ref{eq:8211rs-decomposition})
for $k=2.4$ concerns elementary
functions, for them the chain of transformations~\ref{eq:8211rs-decomposition}
is also working.

As for the Belyi function from Figure~\ref{fig:81111i}, the higher order
transformations of the Gauss hypergeometric function associated with it
cannot be presented (in general) as a composition of transformations of the
lower order, because the first term of the
composition~(\ref{eq:9111-decomposition}) does not have $R_3$-type.
In terms of the $\theta$-triples the corresponding seed transformations
for the Gauss function are as follows:
$$
\left(\frac12,\frac13,\frac{k}9\right)\leftarrow
\left(\frac{k}9,\frac{k}9,\frac{k}9-\epsilon_k\right),\quad
k=1,2,3,4,\quad\epsilon_0=\epsilon_2=1-\epsilon_1=1-\epsilon_3=0.
$$
Because of the composition character both Belyi functions discussed here are
easy to construct. We present here the function~(\ref{eq:9111-decomposition}),
because as follows from the above discussion it is a useful function in the
theory of the Gauss hypergeometric functions:
$$
\lambda_1=\frac{64(\lambda^3-1)}{\lambda^3(9\lambda^3-8)^3},\qquad
\lambda_1=-\frac{192i\sqrt{3}\lambda(\lambda-1)
(\lambda-\frac{1-i\sqrt{3}}2)^9}{(\lambda-\frac{1+i\sqrt{3}}2)^3
(\lambda^3-\frac{3(1+17i\sqrt{3})}2\lambda^2-\frac{3(1-17i\sqrt{3})}2\lambda+1)^3}
$$
The first formula shows how this function is constructed via the composition
given in Equation~(\ref{eq:9111-decomposition}) and the second one -
its much more sophisticated form after a normalization suitable for the
construction of the higher order transformation of the Gauss hypergeometric
function.}\end{remark}

Turning back to the discussion of Figure~\ref{fig:81111i}
it is interesting to notice that although  Equation~(\ref{eq:composition})
formally still holds for the R-types the deformation function is
indecomposibale, see Remark~\ref{rem:81111i} below.
A direct analysis of the system of algebraic equations defining this
rational covering (see \cite{HGBAA}, Remark 2.1) reveals another solution:
\begin{align}
 \label{eq:81111i}
\lambda_1&=\frac{27(s^2+s+is-i)^8}{8s^2(s^4+1)^3}
\frac{\lambda(\lambda-1)(\lambda-t)(\lambda-a)^8}
{(\lambda^4+c_3\lambda^3+c_2\lambda^2+c_1\lambda+c_0)^3},\\
\label{eq:a81111i}
a&=\frac{(1+i)(s^2+1)(s^2+is-s+i)(s^2+2s-1)^3}{8s(s^2+s+is-i)(s^2+i)^3},
\end{align}
\begin{align*}
c_3&=\frac{(s^2+2s-1)}{2s^2(s^2-i)^2(s^2+i)^3}
(s^{12}+(5+7i)s^{11}+(31i-13)s^{10}+(5i-29)s^9+(32+17i)s^8\\
&+(62i-2)s^7+(28-28i)s^6+(62-2i)s^5-(17+32i)s^4+(5-29i)s^3+(13i-31)s^2\\
&+(7+5i)s-i),\\
c_2&=\frac{(s^2+2s-1)^4}{64s^4(s^2-i)^2(s^2+i)^6}
(s^{16}+8is^{15}-68s^{14}-296is^{13}+252s^{12}-184is^{11}+420s^{10}\\
&+472is^9+454s^8+472is^7+420s^6-184is^5+252s^4-296is^3-68s^2+8is+1),
\end{align*}
\begin{align*}
c_1&=-\frac{(s^2+1)^2(s^2+2s-1)^7}{64s^4(s^2-i)^2(s^2+i)^9}
(s^{12}+(7i-5)s^{11}-(13+31i)s^{10}+(29+5i)s^9\\
&+(32-17i)s^8+(2+62i)s^7+(28+28i)s^6-(62+2i)s^5+(32i-17)s^4-(5+29i)s^3\\
&-(31+13i)s^2+(5i-7)s+i),\\
c_0&=\frac{(s^2-2s-1)^2(s^2+1)^4(s^2+2s-1)^{10}}{1024s^4(s^2-i)^2(s^2+i)^{10}},
\end{align*}
The solution obtained from the above function via Theorem 2.1 \cite{HGBAA}
reads,
\begin{align}
 \label{eq:t81111i}
&t=-\frac{(s^2+1)^2(s^2+2s-1)^3(s^2-2s-1)^3}{32s^2(s^2+i)^3(s^2-i)^3},\\
\label{eq:y81111i}
&y(t)=\frac{(1+i)(s^2+s-is+i)(s^2-2s-1)(s^2+1)(s^2+2s-1)^2}
{8s(s^2+i)(s^2-i)^2(s^2-s-is-i))},
\end{align}
It solves Equation~(\ref{eq:P6}) for the following $\theta$-tuple
$$
\theta_0=\theta_1=\theta_t=\frac18,\quad\theta_\infty=\frac78.
$$
\begin{proposition}
 \label{prop:238-RS-81111}
With each rational function~{\rm(\ref{eq:lambda81111})} and 
{\rm(\ref{eq:81111i})} associated four seed $RS$-trans\-for\-mations:\\
$
RS_4^2\left(\!\!
\begin{tabular}{c|c|c}
$k/8$&$1/2$&$1/3$\\
$8+\underbrace{1+\!\ldots\!+1}_4$&$\underbrace{2+\!\ldots\!+2}_6$&
$\underbrace{3+\!\ldots\!+3}_4$
\end{tabular}\!\!
\right)
$
for $k=1,2,3$, and $4$.
Each of these transformations produces one algebraic genus-$0$
solution of Equation~{\rm(\ref{eq:P6})} for the following sets of the 
$\theta$-parameters:
\begin{align*}
k&=1,\qquad&\theta_0&=\theta_1=\theta_t=1-\theta_\infty=\frac18,\\
k&=2,\qquad&\theta_0&=\theta_1=\theta_t=\theta_\infty=\frac14,\\
k&=3,\qquad&\theta_0&=\theta_1=\theta_t=1-\theta_\infty=\frac38,\;\;\text{and}\;\;\\
k&=4,\qquad&\theta_0&=\theta_1=\theta_t=\theta_\infty=\frac12.
\end{align*}
The solutions for $k=1$ corresponding to the functions{\rm(\ref{eq:lambda81111})} and 
{\rm(\ref{eq:81111i})} are given by Equations~{\rm(\ref{eq:t-y-2})} and
{\rm(\ref{eq:t81111i}), (\ref{eq:y81111i})}, respectively. 
\end{proposition}
\begin{remark}
 \label{rem:81111i}{\rm
All solutions of Equation~(\ref{eq:P6}) which can be produced via
Proposition~\ref{prop:238-RS-81111} with the help of the function
(\ref{eq:lambda81111}) are rational functions of $\sqrt{t}$ and $\sqrt{t-1}$.

The situation with the solutions that are generated via the 
function~(\ref{eq:81111i}) is more interesting:
most probably, these solutions (it is not checked yet) coincides with the solution 
that can be obtained via a successive compositions of
$RS_4^2\left(\!\!
\begin{tabular}{c|c|c}
1/4&1/2&1/3\\
4+1+1&2+2+1+1&3+3
\end{tabular}\!\!
\right)$
(See, Subsection 3.3, Example 5 ($CW$-split) of
\cite{HGBAA}) the Okamoto transformation in the sense of Appendix of
\cite{KK}, and one of the quadratic transformations for Equation~(\ref{eq:P6})
at the last page of Appendix of \cite{HGBAA}), such transformations also
can be described as the simple quadratic Belyi functions. Appearance of the
Okamoto transformation in this
composition makes impossible to lift it on the level of rational coverings.

The last statement can be easily observed in this particular case. Because,
In case we suppose that these two examples are related with some
$RS$-transformation, then it would mean that two hypergeometric functions
with the $\theta$-triples $(1/2,1/3,1/8)$, from the above example, and
$(1/2,1/3,1/4$, from Subsection 3.3, Example 5 ($CW$-split) of
\cite{HGBAA}), are related with some algebraic transformation, which is
impossible, because the first function does not belong to the Schwarz
list~\cite{SCH} while the second is in its fourth row.

This, possibly, explains the appearance of $i$ in the parametrization
(\ref{eq:81111i}). Therefore, although the function~(\ref{eq:81111i})
defines an algebraic solution which (most probably) can be obtained from
the already known simplier ones by the certain transformations, the explicit
formula for the covering (\ref{eq:81111i}) has an independent value.
In particular, if we are interesting not only in the solutions of the sixth 
Painlev\'e equation but also in the solutions of the associated linear ODE 
the function~(\ref{eq:81111i}) gives us an additional opportunity to provide 
an explicit construction of the latter function in terms of the hypergeometric ones.}
\end{remark}

\end{document}